# Riemann's Zeta Function.
# An Attempt to fathom ζ(3).

Renaat Van Malderen

November 2012

*Abstract*

Already in 1734 Euler found a short explicit formula for the value of Riemann's zeta function ζ(s) when the argument s equals a positive integer 2n where n=1,2,3…No such formula exists for odd-positive integer arguments of Zeta. The present paper discusses in particular the case of ζ(3). A formula for ζ(3) is obtained which in addition to a number of well known constants includes a rapidly converging infinite series, of which each term contains rational numbers and an even power of π. An attempt to convert this series into a finite number of terms containing commonly known constants is met with only partial success. The general case for ζ(2n+1) is also worked out.

## 1. Introduction.

In a previous paper by the author [1] a formula was given for Riemann's Zeta Function:

$$\zeta(s) = 1 + \frac{1}{s-1}\left(\frac{2}{3}\right)^{s-1} - \sum_{k=1}^{\infty} \frac{\prod_{r=0}^{2k-1}(s+r)[\zeta(s+2k)-1]}{(2k+1)!\, 2^{2k}} \qquad (1)$$

(1) Is absolutely convergent for all complex s except for s=1 where a simple pole occurs. Using (1) to compute actual values of ζ(s) requires knowledge of the terms ζ(s+2k). Nevertheless (1) allows to extend ζ(s) into Re(s)≤1. For s=0 or a negative integer the series (1) reduces to a finite number of terms. The value of ζ(s) for these s values is well known and equals:

$$\zeta(-n) = \frac{(-1)^n B_{n+1}}{(n+1)} \qquad (2)$$

$B_{n+1}$ is the Bernoulli number of index n+1 and of course for s equal a negative even-integer the value of ζ(-n) equals zero; i.e. for n=2m, ζ(-2m)=0. For even positive integers 2m the value of Zeta is known as well (Euler):

$$\zeta(2m) = \frac{|B_{2m}|(2\pi)^{2m}}{2\,(2m)!} \qquad (3)$$

For odd positive integers however a formula similar to (3) is not known. Expressing ζ(3) in terms of (1) results in a series containing ζ(5), ζ(7), … . Below it is shown how to express ζ(3) as a rapidly converging series in which only even, i.e. ζ(2m) values occur and as given in (3), the terms of this series include only rational numbers and even powers of π. An attempt to reduce the mentioned series to a finite expression in terms of commonly known constants such as π, e, and logarithms only met with partial success.

## 2. The Functional Equation.

Riemann's famous functional equation

$$\zeta(1-s) = \pi^{-s} 2^{1-s} \Gamma(s) \zeta(s) \cos\left(\frac{\pi s}{2}\right) \quad (4)$$

expresses ζ(1-s) in terms of ζ(s) and is the key to what follows: Even negative integer values of s translate into odd positive integer values for 1-s. So, if we put s=-2 into (4) we get:

$$\zeta(3) = 8\pi^2 \Gamma(-2) \zeta(-2) \cos(-\pi) = -8\pi^2 \Gamma(-2) \zeta(-2) \quad (5)$$

$\Gamma(s)$ has simple poles for s=0 and all negative integers and ζ(-2)=0. Therefore $\Gamma(-2) \zeta(-2)$ is undefined but we will determine its limiting value by putting s=-2+ε with ε approaching zero. $\Gamma(-2+\varepsilon)$ then keeps getting larger whereas ζ(-2+ε) goes to zero. Nevertheless ζ(3) is well defined. So the product ϕ(-2+ε)=$\Gamma(-2+\varepsilon)$ζ(-2+ε) must approach a finite non zero limit. For convenience in what follows we will often use:

Z(s)= ζ(s)-1   (6).

## 3. The equation for ζ(3).

Expression for ζ(-2+ε).

Plugging s= -2+ε into (1) yields:

$$\zeta(-2+\varepsilon) = 1 + T_0 + T_1 + \sum_{k=2}^{\infty} T_k \quad (7)$$

in which $T_0 = \frac{1}{-3+\varepsilon}\left(\frac{2}{3}\right)^{-3+\varepsilon}$ and neglecting terms in $\varepsilon^2$ and higher:

$$T_0 = -\frac{9}{8}\left[1 + \varepsilon(\ln\left(\frac{2}{3}\right) + \frac{1}{3})\right] \quad (8)$$

For k=1: $T_1$= -(-2+ε)(-1+ε)Z(0+ε)/(3! $2^2$)

Ignoring again terms in $\varepsilon^2$ and higher:

Z(0+ε)= ζ(0+ε)-1=ζ(0)+εζ'(0)-1

As well known: $\zeta(0) = -1/2$ and $\zeta'(0) = \frac{-\ln(2\pi)}{2}$

This yields: $= T_1 = \frac{1}{8} + \varepsilon\left(\frac{\ln(2\pi)}{24} - \frac{3}{16}\right)$ (9).

In the terms $T_k$ with k=2,3,... with the product $\prod_{r=0}^{2k-1}(s+r) = \prod_{r=0}^{2k-1}(-2+r+\varepsilon)$, the factor corresponding to r=2 reduces to ε only. As a result in all other factors, including $Z(-2+2k+\varepsilon)$, ε may be ignored since we don't consider $\varepsilon^2$ and higher.

So:

$$T_k = \frac{(-2)(-1)\varepsilon(1)(2)\ldots(2k-3)Z(2(k-1))}{(2k+1)!\,2^{2k}}$$

$$T_k = \frac{2\varepsilon(2k-3)!\,Z(2(k-1))}{(2k+1)!\,2^{2k}}$$

Switching from k to l=k-1, l=1,2,3,....:

$$\sum_{k=2}^{\infty} T_k = \sum_{l=1}^{\infty} T_l = \sum_{l=1}^{\infty} \frac{2\varepsilon(2l-1)!\,Z(2l)}{(2l+3)!\,2^{2(l+1)}} = 2\varepsilon \sum_{l=1}^{\infty} C_l Z(2l) \quad (10)$$

with

$$C_l = \frac{1}{16}\frac{1}{2^{2l}}\frac{1}{l(l+1)(2l+1)(2l+3)} \quad (11)$$

Adding up $T_0$, $T_1$ and $\sum T_l$:

$$\zeta(-2+\varepsilon) = \frac{9\varepsilon}{8}\left[\ln\left(\frac{3}{2}\right) - \frac{1}{2} + \frac{\ln(2\pi)}{27}\right] - 2\varepsilon \sum_{l=1}^{\infty} C_l Z(2l) \quad (12)$$

Expression for $\Gamma(-2+\varepsilon)$

According to the Gamma function's fundamental property: $\Gamma(s+3) = (s+2)(s+1)s\,\Gamma(s)$
For s= -2+ε:

$$\Gamma(-2+\varepsilon) = \frac{\Gamma(1+\varepsilon)}{\varepsilon(\varepsilon-1)(\varepsilon-2)} = \frac{\Gamma(1) + \varepsilon\Gamma'(1)}{2\varepsilon - \varepsilon^2 + \cdots}$$

$$\Gamma(-2+\varepsilon) = \frac{\Gamma(1)}{2\varepsilon} + \frac{\Gamma'(1)}{2} \quad (13)$$

The product ϕ(-2+ε)

Multiplying (12) and (13) and letting ε go to zero:

$$\phi(-2) = \frac{1}{16}\left[\frac{\ln(2\pi)}{3} + 9\ln\left(\frac{3}{2}\right) - \frac{9}{2}\right] - \sum_{l=1}^{\infty} C_l Z(2l) \qquad (14)$$

As is obvious, the term **Γ'(1)/2** in (13) is inconsequential.

So finally: $\zeta(3) = -8\pi^2 \phi(-2)$

$$\zeta(3) = -8\pi^2 \left[\frac{1}{16}\left(\frac{\ln(2\pi)}{3} + 9\ln\left(\frac{3}{2}\right) - \frac{9}{2}\right) - \sum_{l=1}^{\infty} C_l Z(2l)\right] \qquad (15)$$

### 4. Numerical considerations.

In line with (3):

$$Z(2l) = \frac{|B_{2l}|(2\pi)^{2l}}{2\,(2l)!} - 1 \qquad (16)$$

So the formula for Z(2l) contains only rational numbers and an even power of π. In addition ζ(3) as given by (15) contains also further rationals as well as the constants ln(2), ln(3) and ln(π). The even powers of π however occur as part of an infinite series. This series, i.e. $\sum_{l=1}^{\infty} C_l Z(2l)$ converges rapidly as shown numerically below:

$$\text{Let A} = \frac{1}{16}\left[\frac{\ln(2\pi)}{3} + 9\ln\left(\frac{3}{2}\right) - \frac{9}{2}\right] = -0.01488677114\ldots \quad (17)$$

Considering just the initial five terms (l=1,2,3,4,5) with an accuracy of eleven decimals past the point:

$C_1 Z(2)\ \ = 0.00033590316$
$C_2 Z(4)\ \ = 0.00000153131$
$C_3 Z(6)\ \ = 0.00000002240$
$C_4 Z(8)\ \ = 0.00000000050$
$C_5 Z(10) = 0.00000000001$

$$\text{Let B} = -\sum_{l=1}^{5} C_l Z(l) = -0.00033745738 \qquad (18)$$

A+B = -0.01522422852…

$\zeta(3) = -8\pi^2(A+B) = 1.202056902\ldots$

The more precise value equals $1.202056903\ldots$ [2].
So the difference ($\sim 10^{-9}$) occurs at the 9th decimal after the point.

## 5. Error Quantification

In this section we have a closer look at the obtainable precision in calculating $\zeta(3)$ as given by (15). A first obvious point to be made when it comes to precision in numerical calculations is that apart from necessary truncation of infinite series, also constants such as $\pi$, ln 2, ln3, will enter the calculations with limited precision, dependent on the computational tools used. Below we concentrate on the error resulting from truncating the series $\sum C_l Z(2l)$ to a given maximum $\sum_{l=1}^{m} C_l Z(2l)$. This requires to put a bound on $Z(2l)$.

As well known:

$$Z(2l) = \frac{1}{\Gamma(2l)} \int_0^\infty \frac{\lambda^{2l-1} e^{-\lambda} d\lambda}{1 - e^{-\lambda}} - 1 = \frac{1}{\Gamma(2l)} \int_0^\infty \frac{\lambda^{2l-1} e^{-2\lambda} d\lambda}{1 - e^{-\lambda}}$$

Now:

$$\frac{\lambda}{1 - e^{-\lambda}} \leq \lambda + 1 \text{ for } 0 \leq \lambda < \infty \quad \text{or} \quad \frac{1}{1 - e^{-\lambda}} \leq 1 + \frac{1}{\lambda}$$

So:

$$Z(2l) < \frac{1}{2^{2l}} \left(1 + \frac{2}{2l - 1}\right) = F(l) \quad (19)$$

For large l, Z(l) approaches $2^{-2l}$. Table -1 shows how $Z(2l)$, $F(l)$ and $2^{-2l}$ compare.

table -1

| l | Z(2l) | F(l) | $2^{-2l}$ |
|---|---|---|---|
| 1 | 0.645 | 0.75 | 0.25 |
| 2 | 8.2E(-2) | 10.4E(-2) | 6.25E(-2) |
| 3 | 1.7E(-2) | 2.2E(-2) | 1.56E(-2) |
| 4 | 4E(-3) | 5E(-3) | 3.9E(-3) |
| 5 | 0.99E(-3) | 1.2E(-3) | 0.97E(-3) |
| 10 | 9.53E(-7) | 10E(-7) | 9.53E(-7) |
| 30 | 9.31E(-10) | 9.96E(-10) | 9.31E(-10) |
| 42 | 2.27E(-13) | 2.38E(-13) | 2.27E(-13) |

Truncation Error

(11) implies $C_l < \frac{1}{2^{2l+6}l^4}$. Combined with (19) provides the Absolute truncation Error AE:

$$AE = \sum_{l=m+1}^{\infty} C_l Z(2l) < \left(1 + \frac{2}{2m+1}\right) \sum_{l=m+1}^{\infty} \frac{1}{2^{4l+6} \, l^4} \quad (20)$$

and the Relative truncation Error $RE = \frac{AE}{\sum_{l=m+1}^{\infty} C_l Z(2l)}$ (21)

with $\sum_{l=m+1}^{\infty} C_l Z(2l) \cong 3.38E(-4)$. By limiting the summation in (20) to l=n we obtain table -2:

table -2

| m | n | AE | RE |
|---|---|---|---|
| 4 | 10 | 2.5E(-11) | 7.4E(-9) |
| 5 | 11 | 7.4E(-13) | 2.2E(-9) |
| 6 | 12 | 2.5E(-14) | 7.4E(-11) |

Rate of convergence

Comparing the ratio of successive terms, i.e.:

$$R(l) = \frac{C_l Z(2l)}{C_{l+1} Z(2(l+1))} \quad (22)$$

gives a good idea of the rate at which the series $\sum C_l Z(2l)$ converges. This is shown in table -3.

table -3

| l | R(l) |
|---|---|
| 1 | 219 |
| 2 | 68 |
| 3 | 44 |
| 4 | 36 |
| 5 | 30 |

For increasing l, R(l) keeps going down but in the limit as l→∞, it still approaches R(l) = 16 as can be seen from (11) and using $Z(2l) \cong \frac{1}{2^{2l}}$.

## 6. Further reduction to finite expressions.

We will now investigate to what extent we might convert the series $\sum C_l Z(2l)$ as given in (15) to a finite number of terms. Partial fraction expansion of (11) yields:

$$C_l = \frac{1}{16} \frac{1}{2^{2l}} \left[\frac{1}{3l} - \frac{2}{2l+1} + \frac{1}{l+1} - \frac{2}{3(2l+3)}\right] \quad (23)$$

As a preliminary we introduce two formulas which we will need below:

$$\sum_{l=1}^{\infty} \frac{\zeta(2l)}{(2l+1)4^l} = \frac{1-\ln 2}{2} \quad (24)$$

$$\sum_{l=1}^{\infty} \frac{\zeta(2l)}{l 4^l} = \ln\left(\frac{\pi}{2}\right) \quad (25)$$

Both (24) and (25) may be obtained by expanding $\sum_{l=1}^{\infty}\left[1 + l \ln\left(\frac{1-\frac{l}{2}}{1+\frac{l}{2}}\right)\right]$ and $-\frac{1}{2}\sum_{l=1}^{\infty}\left(1 - \frac{1}{4l^2}\right)$ into power series and by obtaining the constants $\frac{(1-\ln 2)}{2}$ and $\ln\left(\frac{\pi}{2}\right)$ as by-products of proving Stirling's formula.

In line with (23) we split $\sum C_l Z(2l)$ into 4 parts:

$$SA = \frac{1}{16}\frac{1}{3}\sum_{l=1}^{\infty}\frac{1}{l 4^l}(\zeta(2l) - 1) \quad (26)$$

$$SB = -\frac{1}{8}\sum_{l=1}^{\infty}\frac{1}{(2l+1)4^l}(\zeta(2l) - 1) \quad (27)$$

$$SC = \frac{1}{16}\sum_{l=1}^{\infty}\frac{1}{(l+1)4^l}(\zeta(2l) - 1) \quad (28)$$

$$SD = -\frac{1}{8}\frac{1}{3}\sum_{l=1}^{\infty}\frac{1}{(2l+3)4^l}(\zeta(2l) - 1) \quad (29)$$

Reduction of SA

Using (25) $\sum_{l=1}^{\infty}\frac{\zeta(2l)}{l 4^l} = \ln\left(\frac{\pi}{2}\right)$

Using $\ln\left(\frac{1}{1-x}\right) = \sum_{l=1}^{\infty}\frac{x^l}{l}$ with $x = \frac{1}{4}$ we have :

$\sum_{l=1}^{\infty}\frac{1}{l 4^l} = \ln\left(\frac{4}{3}\right)$.

Combining these two results :

$$SA = \frac{1}{48}\left(\ln\left(\frac{\pi}{2}\right) - \ln\left(\frac{4}{3}\right)\right) \quad (30)$$

Reduction of SB

Using (24): $-\frac{1}{8}\sum_{l=1}^{\infty}\frac{\zeta(2l)}{(2l+1)4^l} = -\frac{1}{16}(1 - ln2)$

Using $\frac{1}{2x}ln\left(\frac{1+x}{1-x}\right) - 1 = \sum_{l=1}^{\infty}\frac{x^{2l}}{(2l+1)}$ with $x = \frac{1}{2}$ we obtain:

$$\frac{1}{8}\sum_{l=1}^{\infty}\frac{1}{(2l+1)4^l} = \frac{1}{8}(ln3 - 1)$$

Combining both results :

$$SB = \frac{1}{8}\left[ln(3\sqrt{2}) - \frac{3}{2}\right] \quad (31)$$

The case of SC

Using $-\left[\frac{ln(1-x)}{x} + 1\right] = \sum_{l=1}^{\infty}\frac{x^l}{l+1}$ and putting $x = \frac{1}{(2n)^2}$ we have:

$$SC = \frac{1}{16}\sum_{l=1}^{\infty}\frac{1}{l+1}\sum_{n=2}^{\infty}\frac{1}{(2n)^{2l}}$$
$$= \frac{1}{16}\sum_{n=2}^{\infty}\sum_{l=1}^{\infty}\frac{1}{l+1}\frac{1}{(2n)^{2l}} = -\frac{1}{16}\sum_{n=2}^{\infty}ln\left[e\left(1-\frac{1}{4n^2}\right)^{4n^2}\right]$$

Or written in product form:

$$SC = -\frac{1}{16}ln\left[\prod_{n=2}^{\infty}e\left(1-\frac{1}{4n^2}\right)^{4n^2}\right] \quad (32)$$

Formula (32) may be of theoretical interest but it is not suitable for numerical evaluation. On the other hand Z(2l) is easily obtained to any desired degree of accuracy (see (3) and (6)). Therefore direct evaluation of (28) yields much better accuracy with far less effort. E.g; using l=1,…, 5, (28) yields SC=0.005150.
(32) only comes close to this result and yields SC=0.005145… after multiplying 1499 terms (n=2,…,1500)! This is due to the slow convergence of the basic series $\sum\frac{1}{n^2}$.

The case of SD

Using $\frac{1}{2x^3}ln\left(\frac{1+x}{1-x}\right) - \frac{1}{x^2} - \frac{1}{3} = \sum_{l=1}^{\infty}\frac{x^{2l}}{2l+3}$ and putting $x = \frac{1}{2n}$ we have:

$$4n^3 \ln\left(\frac{1+\frac{1}{2}n}{1-\frac{1}{2}n}\right) - 4n^2 - \frac{1}{3} = \sum_{l=1}^{\infty} \frac{1}{(2l+3)(2n)^{2l}}$$

$$SD = -\frac{1}{6}\sum_{n=2}^{\infty}\left(n^3\ln\left(\frac{1+\frac{1}{2}n}{1-\frac{1}{2}n}\right) - n^2 - \frac{1}{12}\right)$$

Or written as an infinite product:

$$SD = -\frac{1}{6}\ln\left[\prod_{n=2}^{\infty}\left[e^{-(n^2+\frac{1}{12})}\left(\frac{1+\frac{1}{2}n}{1-\frac{1}{2}n}\right)^{n^3}\right]\right] \quad (33)$$

Similar to the comment for SC, also for SD, formula (29) converges rapidly while (33) converges very slowly. Checking (SA+SB+SC+SD) against $\sum_{l=1}^{5} C_l Z(2l)$ as given in (18):
Using (30), (31), (28) and (29) we obtain:

SA =  0.003414596…
SB = -0.006851765…
SC =  0.005150182…
SD = -0.001375556…
______________________

Σ = 0.000337457… which equals $\sum_{l=1}^{5} C_l Z(2l)$ as obtained in (18) up to the 9th decimal after the point.
Incorporating the expressions for SA and SB, as given by (30) and (31) into (15) results in:

$$\zeta(3) = -8\pi^2\left(\frac{5}{12}\ln 3 - \frac{13}{24}\ln 2 - \frac{3}{32} - SC - SD\right) \quad (34)$$

with SC and SD as given by either (28), (29) or (32), (33).

### 7. The general case of ζ(2n+1)

In this section we will generalize the approach taken for the ζ(3) case. As a preliminary we consider the curve ζ(σ) for σ real and negative.
For odd negative integers:

$$\zeta(-2n+1) = \frac{(-1)^n(2n-1)!\,\zeta(2n)2}{(2\pi)^{2n}} \quad n = 1,2,\ldots \quad (35)$$

For even negative integers ζ(-2n)=0.
ζ(σ) in the range -∞<σ<0 is a wavy line with zeros at σ=-2n and the bulges in between the zeros having an amplitude at σ=-2n+1 (n=2,3,…) given by (35), alternating in the positive and negative direction.

The minimum absolute height of the bulges occurs for n=3 and 4. For larger n values $|\zeta(-2n+1)|$ keeps growing. This means that at the zero crossings (i.e. σ=-2n) the curve ζ(σ) will get increasingly steeper.

Equation (9) for the case of ζ(3) involved the value of ζ'(0). As we will find out below, for the general case of ζ(2n+1), we will need to know in addition to ζ'(0), also ζ'(-2l) for l=1,2,…, up to l=n-1.

## 8. The Derivative ζ'(-2l)

Using again the functional equation

$$\zeta(1-s) = \pi^{-s} 2^{1-s} \,\Gamma(s)\zeta(s)\cos\left(\frac{\pi s}{2}\right) \qquad (36)$$

and inserting s=-2l we again get the undefined expression $\Gamma(-2l)\zeta(-2l)$.

Similar to the case for ζ(3) we determine its limiting value. Using similar symbolism as before, now with s=-2l+ε, it is easily established that:

$$\Gamma(-2l+\varepsilon) = \frac{1}{\varepsilon(2l)!} \qquad (37)$$

Since we are after ζ'(-2l) we simply put

$$\zeta(-2l+\varepsilon) = 0 + \varepsilon\zeta'(-2l) \qquad (38)$$

Further, for s=-2l+ε:

$$\cos\left(\frac{\pi s}{2}\right) = \cos\left(-\pi l + \frac{\pi\varepsilon}{2}\right) = (-1)^l \qquad (39)$$

since here we can ignore the ε term.

$$\lim_{\varepsilon\to 0} \Gamma(-2l+\varepsilon)\zeta(-2l+\varepsilon) = \frac{1}{\varepsilon(2l)!}\,\varepsilon\zeta'(-2l) = \frac{\zeta'(-2l)}{(2l)!} \qquad (40)$$

Taking account of (39) and (40), (36) yields:

$$\zeta'(-2l) = \frac{(2l)!\,(-1)^l \zeta(2l+1)}{2\,(2\pi)^{2l}} \qquad (41)$$

In line with the observation about the growing bulges in the curve ζ(σ), also $|\zeta'(-2l)|$ grows increasingly steeper with l going up.

### 9. The Equation for ζ(2n+1)

Neglecting higher order terms, for s=-2n:

$$\zeta(-2n + \varepsilon) = \varepsilon \zeta'(-2n) \qquad (42)$$

This is the term to be used in the functional equation in order to extract from it an expression for $\zeta(2n + 1)$.

We determine ζ'(s) based on equation (1):

$$\zeta'(s) = \frac{d}{ds}\left\{1 + \frac{1}{s-1}\left(\frac{2}{3}\right)^{s-1} - \sum_{k=1}^{\infty} \frac{\prod_{r=0}^{2k-1}(s+r)\, Z(s+2k)}{(2k+1)!\, 2^{2k}}\right\} \qquad (43)$$

We break up (43) into four terms $T_0$, $T_1$, $T_2$, $T_3$ and evaluate each one separately at s=-2n.

Term $T_0$: This one is straightforward: For s=-2n we have

$$T_0 = \frac{d}{ds}\left[1 + \frac{1}{s-1}\left(\frac{2}{3}\right)^{s-1}\right] = \frac{1}{2n+1}\left(\frac{3}{2}\right)^{2n+1}\left[\ln\left(\frac{3}{2}\right) - \frac{1}{2n+1}\right] \qquad (44)$$

Terms $T_1$, $T_2$, $T_3$: The summation over k in (43) is split up in 3 parts, k=1,…., (n-1) for $T_1$, k=n for $T_2$, k=(n+1), …., ∞ for $T_3$. Table -4 shows the relevant expressions to work out T1, T2, T3 evaluated at s=-2n. The term A(n,k) as used in table -4 is defined as:

$$A(n, k) = \sum_{r=0}^{2k-1} \frac{1}{2n - r}$$

Table -4

| k | 1,…., n-1 | n | n+1, …., ∞<br>k=n+l<br>l=1,…..,∞ |
|---|---|---|---|
| Z(-2n+2k) | ζ (-2n+2k)-1 = -1 | $\zeta(0) - 1 = -\frac{3}{2}$ | Z(2l) |
| $\frac{d}{ds}\left[\prod_{r=0}^{2k-1}(s+r)\right]$ | $\frac{-(2n)!\, A(n,k)}{(2n-2k)!}$ | -(2n)!A(n,n) | (2n)!(2l-1)! |
| $\prod_{r=0}^{2k-1}(s+r)$ | $\frac{(2n)!}{(2n-2k)!}$ | (2n)! | 0 |
| $\frac{d}{ds}[Z(s+2k)]$ | ζ'(-2n+2k) | $\zeta'(0) = -\frac{\ln(2\pi)}{2}$ | Irrelevant |

From table -4:

$$T_1 = \sum_{k=1}^{n-1} (2n)! \frac{[A(n,k) + \zeta'(-2n+2k)]}{(2n-2k)!(2k+1)!\, 2^{2k}} \qquad (45)$$

$$T_2 = (2n)! \left[ \frac{3A(n,n) - \ln(2\pi)}{(2n+1)!\, 2^{2n+1}} \right] \qquad (46)$$

$$T_3 = \frac{(2n)!}{2^{2n}} \sum_{l=1}^{\infty} \frac{Z(2l)(2l-1)!}{(2n+2l+1)!\, 2^{2l}} = \frac{(2n)!}{2^{2n}} \sum_{l=1}^{\infty} \frac{Z(2l)}{2^{2l} \prod_{r=0}^{2n+1}(2l+r)} \qquad (47)$$

Taking into account (37), (42):

$$\phi(-2n) = \lim_{\varepsilon \to 0} \Gamma(-2n+\varepsilon)\zeta(-2n+\varepsilon) = \frac{T_0 - T_1 - T_2 - T_3}{(2n)!} \qquad (48)$$

Evaluating the terms $\pi^{-s} 2^{1-s}$ and $\cos\left(\frac{\pi s}{2}\right)$ in (36) at s=-2n and inserting the result together with (48) into (36) we finally get:

$$\zeta(2n+1) = (-1)^n \pi^{2n} 2^{1+2n} \left\{ \frac{1}{(2n+1)!} \left(\frac{3}{2}\right)^{2n+1} \left[\ln\left(\frac{3}{2}\right) - \frac{1}{2n+1}\right] \right.$$

$$- \sum_{k=1}^{n-1} \frac{[A(n,k) + \zeta'(-2n+2k)]}{(2n-2k)!(2k+1)!\, 2^{2k}} - \frac{[3A(n,n) - \ln(2\pi)]}{(2n+1)!\, 2^{2n+1}}$$

$$\left. - \frac{1}{2^{2n}} \sum_{l=1}^{\infty} \frac{Z(2l)}{2^{2l} \prod_{r=0}^{2n+1}(2l+r)} \right\} \qquad (49)$$

Remarks concerning equation (49)

a) Is (49) a "simple" formula? Not really, but that's the way it is. (49) is recursive in the sense that the expression for $\zeta(2n+1)$ contains the terms $\zeta'(-2n+2k)$ with k=1,....,2(n-1), which as given by (41) require the knowledge of $\zeta(2n-1)$, ..., $\zeta(3)$. This is a complication which does not occur for $\zeta(3)$. Notice that also formula (1) is recursive in the same sense.

b) For n=1, i.e; the case of $\zeta(3)$, the second term in (49) drops out (meaningless).

c) Convergence Rate.
The convergence rate of the infinite series in (49), as defined in (22) for $\zeta(3)$ now involves 2 parameters, i.e. n and l. We call it R(n,l). Although for large l, R(n,l) approaches the same value 16 as for $\zeta(3)$, the initial values (rounded off) increase with increasing n (see table -5).

$$R(n,l) = \frac{4Z(2l)(l+n+1)(2l+2n+3)}{Z(2l+2)l(2l+1)} \qquad (50)$$

table -5

|   | n |   |   |   |
|---|---|---|---|---|
| l | 1 | 2 | 3 | 4 |
| 1 | 219 | 376 | 575 | 815 |
| 2 | 68 | 104 | 148 | 199 |
| 3 | 44 | 63 | 85 | 110 |
| 4 | 36 | 48 | 62 | 78 |

Numerical Examples.

a) n=2, ζ(2n+1)= ζ(5).
The terms $T_0$, $T_1$, $T_2$, $T_3$, as given in (44), (45), (46), (47) equal:

$T_0 = 0.312050132939$
$T_1 = 0.276442438138$
$T_2 = 0.027575768335$
$T_3 = 0.000048115016$

$T_3$ includes the first 4 terms (l=1,2,3,4).

$T_0$- $T_1$- $T_2$- $T_3$ = 0.00798381145 = ζ'(-4) as calculated using (44), (45), (46), (47).

Checking the above result for ζ'(-4):
Calculating ζ'(-4) based on (41) and using the ζ(5) value from existing tables (e.g. from [2]) shows that the above obtained value is correct up to the 11th decimal past the point.

(49) then yields:

$$\zeta(5) = \frac{(-1)^2 \pi^4 2^5 \zeta'(-4)}{4!} = 1.03692775514 \dots$$

b) n=3, ζ(2n+1)= ζ(7).
A similar calculation as in a. yields:

ζ(7) = 1.008349 … accurate to 6 decimals past the point, using 11 decimals accuracy for obtaining the expression between the large brackets of (49). This drop in accuracy is caused as follows: On the left side of (49) we find ζ(2n+1) which for all practical purposes equals 1. On the right side of (49) the combination $\pi^{2n} 2^{1+2n}$ keeps growing with n. So the remaining factor on the same side must go down. To obtain the same overall accuracy this factor must become more accurate (more significant digits).

## 10. Concluding remarks

In a recent book [3] which first mentions the formula for ζ(2n), it is subsequently stated: "What is not assured is the nature of ζ(2n+1) for n≥1, since neither Euler nor any who have followed him were able to find a workable expression for the number, much less to prove it to be irrational, until 1978 when … R. Apéry made a series of … assertions which combined to a proof that ζ(3) is irrational, although a workable expression for the number remains elusive". Another book [4] states "As yet no simple formula ….. is known for ζ(2n+1) or even for any special case such as ζ(3)….". Or "No one has yet succeeded in obtaining a formula for ζ(2n+1) as simple as for ζ(2n). In fact besides the result of R. Apéry that ζ(3) is irrational, almost nothing is known about the arithmetical structure of ζ(2n+1)" [5].

Terms such as "workable", "simple", and "arithmetical structure" are open to interpretation. Certainly the expression $\sum_{n=1}^{\infty} \frac{1}{n^3}$ as definition of ζ(3) is simple and could be termed as giving insight into its arithmetical structure although it is hardly suitable (workable) for precise calculation due to its slow convergence. In this respect formula (15) might be considered workable because the infinite series included is rapidly convergent and as argued earlier in this paper, the factors Z(2l) are readily computable and as such (15) might arguably be called simple.

Addendum: The derivative ζ'(-2n+1)

The "bulges" occurring in ζ(σ) for σ<0 are not attaining their maximum amplitude precisely for σ=-2n+1. This is demonstrated below. Using the functional equation and inserting s=2n+ε, yields:

$$\zeta(-2n+1-\varepsilon) = \frac{2(2\pi)^{-\varepsilon}}{(2\pi)^{2n}} \Gamma(2n+\varepsilon)\zeta(2n+\varepsilon)\cos\left[\frac{\pi(2n+\varepsilon)}{2}\right]$$

Disregarding higher order terms $\varepsilon^2$, …, and noting

$$\zeta(-2n+1) = \frac{2(-1)^n}{(2\pi)^{2n}} \Gamma(2n)\zeta(2n) \quad \text{we are left with the ε terms:}$$

$$\zeta'(-2n+1) = \frac{2(-1)^{n+1}}{(2\pi)^{2n}} [\Gamma(2n)(\zeta'(2n) - \zeta(2n)\ln(2\pi)) + \Gamma'(2n)\zeta(2n)]$$

Using the Digamma function ψ(2n) = **Γ**'(2n)/**Γ**(2n):

$$\zeta'(-2n+1) = \frac{2(-1)^{n+1}}{(2\pi)^{2n}} \Gamma(2n)[\zeta'(2n) + \zeta(2n)(\psi(2n) - \ln(2\pi))]$$

Formulas for ψ(2n) exist and have been tabulated [2, p.258,272]. ψ(2n) increases with increasing n. Also ζ'(2n)<0 and approaching negative zero with increasing n. Using the ψ(2n) tables: ψ(2n)-ln(2π)<0 for n=1, 2, 3, while for n≥4 it is >0. Since for n=4 we have ζ'(2n)=ζ'(8)=-0.0029019… and ψ(8) -ln(2π)=+0.1777644…, we may conclude: ζ'(-1), ζ'(-3), ζ'(-5) are all negative, while ζ'(-2n+1) for n≥4 are all positive.

## References

[1]     arXiv: 1109.6790v1[math.NT].
[2]     M. Abramowitz and I. Stegun, Handbook of Mathematical Functions (p. 811), Dover 1990.
[3]     Julian Havil, The Irrationals (p. 138), Princeton University Press, 2012.
[4]     Tom.M. Apostol, Introduction to Analytic Number Theory (p. 267), Springer (corrected fifth printing), 1998.
[5]     A. Ivić, The Riemann Zeta Function (p. 4), Dover, 2003.

Renaat Van Malderen
Address: Maxlaan 21, B-2640 Mortsel, Belgium
The author can be reached by email: hans.van.malderen@telenet.be